\documentclass[12pt]{amsart}
\usepackage{amsmath,amssymb,amsthm,upref,graphicx,mathrsfs}

\usepackage{color}
\usepackage[
  colorlinks=true,
  linkcolor=blue,
  citecolor=blue,
  urlcolor=blue]{hyperref}

\numberwithin{equation}{section}

\newtheorem{theorem}{Theorem}[section]
\newtheorem{prop}[theorem]{Proposition}
\newtheorem{lemma}[theorem]{Lemma}
\newtheorem{cor}[theorem]{Corollary}

\theoremstyle{definition}
\newtheorem{definition}[theorem]{Definition}

\newtheorem{remark}[theorem]{Remark}

\newcommand{\be}{\begin{eqnarray*}}
\newcommand{\ee}{\end{eqnarray*}}
\newcommand{\beq}{\begin{equation}}
\newcommand{\eeq}{\end{equation}}

\begin{document}

\title[Multisublinear maximal operator in martingale spaces]
  {Weighted norm inequalities
for multisublinear maximal operator in martingale spaces}

\authors

\author[W. Chen]{Wei Chen}
\address{Wei Chen \\ School of Mathematical Sciences,
Yangzhou University, 225002 Yangzhou, China}
\email{weichen@yzu.edu.cn}

\author[P. D. Liu]{Peide Liu}
\address{Peide Liu\\
School of Mathematics and Statistics, Wuhan University, 430072 Wuhan , China}
\email{pdliu@whu.edu.cn}

\makeatletter
\renewcommand{\@makefntext}[1]{#1}
\makeatother \footnotetext{\noindent
 Supported by the National Natural
Science Foundation of China (Grant Nos. 11101353 and 11071190),
the Natural Science Foundation of Jiangsu Education Committee (Grant No. 11KJB110018) and the Natural Science Foundation of Jiangsu Province (Grant No. BK2012682)}


%
%
\begin{abstract} {Let $v,~\omega_1, ~\omega_2$ be weights and $1<p_1, ~p_2<\infty.$
Suppose that
$\frac{1}{p}=\frac{1}{p_1 }+\frac{1}{p_2 }$ and $(\omega_1, \omega_2)\in RH(p_1, p_2).$
For the multisublinear maximal operator $\mathfrak{M}$ in martingale spaces, we
characterize the weights for which $\mathfrak{M}$ is bounded from
$L^{p_1}(\omega_1)\times L^{p_2}(\omega_2 )$
to $L^{p,\infty}(v)\hbox{ or }L^p(v).$ If $v=\omega_2^{\frac{p}{p_2}}\omega_2^{\frac{p}{p_2}},$
we partially give the bilinear version of one-weight theory.}
\end{abstract}

\maketitle

%
%
 \section{Introduction }
Let $R^n$ be the $n\hbox{-dimensional}$ real Euclidean
space and $f$ a real valued measurable function, the classical
Hardy-Littlewood maximal operator $M,$ the maximal geometric mean
operator $G$ and the minimal operator $\mathrm{m}$ are defined by
$$Mf(x)=\sup\limits_{x\in Q}\frac{1}{|Q|}\int_Q|f(y)|dy,$$
$$G(f)(x)=\sup\limits_{x\in Q} \exp{\frac{1}{|Q|}}\int_Q \log |f(y)|dy$$
 and $$\mathrm{m}f(x)=\inf\limits_{x\in Q}\frac{1}{|Q|}\int_Q|f(y)|dy.$$
 where $Q$ is a non-degenerate cube with its sides parallel to the coordinate
axes and $|Q|$ is the Lebesgue measure of $Q.$

Let $u,~v$ be two weights, i.e. positive measurable functions. As
well known, for $p\geq1,$ \cite{Muckenhoupt B} showed that the
inequality
$$ \lambda^p\int_{\{Mf>\lambda\}}u(x)dx
\leq C \int_{R^n}|f(x)|^pv(x)dx, ~~\lambda>0,~f\in{L^p(v)}
$$
holds if and only if $(u,v)\in A_p,$ i.e., for any cube $Q$ in
$R^n$ with sides parallel to the coordinates
$$
\big(\frac{1}{|Q|}\int_Qu(x)dx\big)
\big(\frac{1}{|Q|}\int_Qv(x)^{-\frac{1}{p-1}}dx\big)^{p-1}<C,~p>1;
$$
$$
\frac{1}{|Q|}\int_Qu(x)dx\leq C \mathop{\hbox{ess inf}}\limits_{Q}v(x)
,~p=1.
$$
Suppose that $u=v$ and $p>1,$ \cite{Muckenhoupt B} also proved that
$$ \int_{R^n}\big(Mf(x)\big)^pv(x)dx
\leq C \int_{R^n}|f(x)|^pv(x)dx, ~\forall f\in{L^p(v)}
$$
holds if and only if $v$ satisfies
\begin{equation}\label{Ap}
\big(\frac{1}{|Q|}\int_Qv(x)dx\big)
\big(\frac{1}{|Q|}\int_Qv(x)^{-\frac{1}{p-1}}dx\big)^{p-1}<C,~\forall Q.
\end{equation}
The crucial step is to show that if $v$ satisfies $A_p$,
then there is an $\varepsilon> 0$ such that $v$
also satisfies $A_{p-\varepsilon}.$
But, the problem of finding all $u$ and $v$ such that
$$ \int_{R^n}\big(Mf(x)\big)^pu(x)dx
\leq C \int_{R^n}|f(x)|^pv(x)dx,~\forall f\in{L^p(v)}
$$
is much hard and complicated. In order to solve the problem,
\cite{Sawyer E T.} established the
testing condition $S_{p,q},$ i.e. for any cube $Q$ in
$R^n$ with sides parallel to the coordinates
$$\Big(\int_{Q}\big(M(\chi_Qv^{1-p'})(x)\big)^qu(x)dx\Big)^{\frac{1}{q}}
\leq C(\int_Qv(x)^{1-p'}dx)^{\frac{1}{p}},~\forall Q$$
where $1<p\leq q<\infty.$
The condition $S_{p,q}$ is a
sufficient and necessary condition such that
the weighted inequality
$$\Big(\int_{R^n}\big(Mf(x)\big)^qu(x)dx\Big)^{\frac{1}{q}} \leq C
\Big(\int_{R^n}|f(x)|^pv(x)dx\Big)^{\frac{1}{p}}, ~\forall f\in{L^p(v)}$$
holds.
In this case, the method of proof is very interesting.
Motivated by \cite{Muckenhoupt B, Sawyer E T.},
the theory of weights developed so
rapidly that it is difficult to give its history a full
account here (see \cite{Garcia Rubio}
and \cite{Cruz-Uribe D. J. M. Martell} for more information).
However, it is possible to give a story of
weighted inequalities for the
different variants of Hardy-littlewood operator.
Let $p\rightarrow\infty$ in $\eqref{Ap},$ it follows that
\begin{equation}\label{A00} \big(\frac{1}{|Q|}\int_Qv(x)dx\big)
\exp\big(\frac{1}{|Q|}\int_Q\log(\frac{1}{v(x)})dx\big)<C,\end{equation}
which is an alternative definition of $A_\infty$ weight (see
\cite{Hruscev S. V.}). It is known that \cite{Sbordone C. and Wik I.}
used $\eqref{A00}$ to characterize the
boundedness of $G$ from $L^1{(v)}$ to $L^1{(v)}.$ In the case of two
weights, \cite{Yin X. Q. and Muckenhoupt B.} gave that
$$\big(\frac{1}{|Q|}\int_Qu(x)dx\big)
\exp\big(\frac{1}{|Q|}\int_Q\log(\frac{1}{v(x)})dx\big)<C,~\forall~
Q\Leftrightarrow\sup\limits_{\parallel
f\parallel_{L^p(v)}=1}\parallel
Gf\parallel_{L^{p,\infty}(u)}<\infty$$ and
$$\int_QG(v^{-1}\chi_Q)(x)u(x)dx\leq C|Q|,
~\forall~ Q\Leftrightarrow\sup\limits_{\parallel
f\parallel_{L^p(v)}=1}\parallel
Gf\parallel_{L^p(u)}<\infty,$$ which generalize the
results of \cite{Wei H. Shi X. L. and Sun Q. Y.}. Recently,
 \cite{Cruz-Uribe D.} (see also the references therein) also studied
the minimal operator and reverse Holder's inequality.
There are still other variants of Hardy-littlewood operator, for example, the
generalized maximal operator and the strong maximal operator
which were considered in
\cite{F. Ruiza, F. Ruiz} and \cite{Kurtz D. S.}, respectively.
Now, the multisublinear maximal function
\begin{equation}\label{multi_maximal_operator}\mathcal{M}(f_1,...,f_m)(x) = \sup\limits_{x\in Q}
\prod\limits_{i=1}\limits^{m}\frac{1}{|Q|}\int_Q|f_i(y_i)|dy_i
\end{equation}
associated with cubes with sides parallel to the coordinate
axes was studied in \cite{Lerner A.K. Ombrosi S.}. They
introduced the multilinear $A_{\overrightarrow{p}}$
condition which is an analogue of the
$A_p$ weight for multiple weights.
The more general case was extensively discussed in
\cite{Grafakos L. Liu L. G. Yang D. C., Grafakos L. Liu L. G. Perez C. Torres R. H.}.

The above operators can be defined in martingale space, and
the weighted inequalities also have their martingale
versions. In fact, all of them have been discussed in
\cite{Zuo H. L. and Liu P. D., Long R. L. and Peng L. Z.,
X. Q. Chang, M. Izumisawa, W. Chen, R. L. Long} (see also the references therein),
except the one for multisublinear maximal function.
In this paper, with stopping times and a kind of reverse Holder's condition,
we discuss weighted inequalities
for multisublinear maximal operator in martingale spaces.
One of our main results is the martingale-variant of $A_{\overrightarrow{p}},$ and
the other is the equivalence of $S_{\overrightarrow{p}}$
and strong weighted inequality in martingale space.
We also discuss the convergence of martingale,
which is partly a bilinear version of the results in \cite{M. Kikuchi}.

The rest of this section consists of the preliminaries for our paper.

Let $(\Omega,\mathcal {F},\mu)$ be a complete probability space and
$(\mathcal {F}_n)_{n\geq0}$ an increasing sequence of
sub$\hbox{-}\sigma\hbox{-}$fields of $\mathcal{F}$ with
$\mathcal{F}=\bigvee\limits_{n\geq0}\mathcal{F}_n.$ A weight
$\omega$ is a random variable with $\omega>0$ and $
E(\omega)<\infty.$ For any $n\geq0$ and $f\in L^1,$ we denote the conditional expectation
with respect to $\mathcal {F}_n$ by $E_n(f),$ $E(f|\mathcal {F}_n)$ or $f_n,$ then $(f_n)_{n\geq0}$ is
an uniformly integral martingale. Suppose that functions $f,~g$ are
integral, the maximal operator
and multisublinear maximal operator are defined by
$$Mf=\sup_{n\geq 0}|E_n(f)|\text{ and }
\mathcal{M}(f,g)=\sup_{n\geq 0}|E_n(f)||E_n(g)|,$$
respectively. Let $B\in \mathcal {F},$
we always denote $\int_\Omega\chi_Bd\mu$ and
$\int_\Omega\chi_B\omega d\mu$ by $|B|$ and $|B|_\omega,$
respectively. For $(\Omega,\mathcal {F},\mu)$ and $(\mathcal
{F}_n)_{n\geq0},$ the family of all stopping times is denoted by
$\mathcal {T}.$ Throughout this paper, $C$ will denote a
constant not necessarily the same at each occurrence.

\thanks{\textbf{Acknowledgement.} This paper was completed while
the first author was at the Faculty of Mathematics of the
University of Seville, Spain. He is very grateful for the hospitality. We also thank Gang Li for many valuable comments on this paper.}
\section{Results and Their Proofs }
\begin{definition}Let $\omega_1, ~\omega_2$ be weights and $1<p_1, ~p_2<\infty.$
Suppose that
$\frac{1}{p}=\frac{1}{p_1 }+\frac{1}{p_2 }$
and $\sigma_i=\omega_i^{-\frac{1}{p_i-1}}\in L^1,$ $i=1,~2.$
We say that the couple of weights $(\omega_1, ~\omega_2)$
satisfies the reverse Holder's condition $RH(p_1,p_2 ),$ if
there exists a positive constant $C$ such that
\be\big(\int_{\{\tau<\infty\}}\sigma_1d\mu\big)^{\frac{p}{p_1}}
\big(\int_{\{\tau<\infty\}}\sigma_2d\mu\big)^{\frac{p}{p_2}}
\leq C\int_{\{\tau<\infty\}}\sigma_1^{\frac{p}{p_1}}
\sigma_2^{\frac{p}{p_2}}d\mu,~\forall \tau\in\mathcal{T}.\ee
\end{definition}

\begin{remark}
In literature there exist many inverse Holder's inequalities of the type
$$\|f\|_p\|g\|_q\leq C\|fg\|,$$
where $\frac{1}{p}+\frac{1}{q}=1,$ $C$ is a constant and the functions $f$ and $g$ are subjected to suitable
restrictions. The suitable restrictions can be found in \cite{Nehari,Zhuang}.
In our paper, we find that the reverse Holder's condition is useful for bilinear weighted theory in martingale context.
\end{remark}

\begin{definition}\label{definition Ap}Let $v,~\omega_1, ~\omega_2$ be weights and $1<p_1, ~p_2<\infty.$
Suppose that
$\frac{1}{p}=\frac{1}{p_1 }+\frac{1}{p_2 }.$
Denote that $\overrightarrow{p}=(p_1,p_2)$
and $\sigma_i=\omega_i^{-\frac{1}{p_i-1}}\in L^1,~i=1,~2.$
We say that the triple of weights $(v,~\omega_1, ~\omega_2)$
satisfies the condition $A_{\overrightarrow{p}},$ if
there exists a positive constant $C$ such that
      \be
\sup\limits_{n\geq0}E_n(v)^{\frac{1}{p}}E_n(\omega_1^{1-p'_1})^{\frac{1}{p'_1}}
      E_n(\omega_2^{1-p'_2})^{\frac{1}{p'_2}}\leq C,
      \ee
where $\frac{1}{p_i}+\frac{1}{p'_i}=1,~i=1,~2.$
\end{definition}

\begin{definition}\label{definition Sp}Let $v,~\omega_1, ~\omega_2$ be weights and $1<p_1, ~p_2<\infty.$
Suppose that
$\frac{1}{p}=\frac{1}{p_1 }+\frac{1}{p_2 }.$ Denote that $\overrightarrow{p}=(p_1,p_2)$
and $\sigma_i=\omega_i^{-\frac{1}{p_i-1}}\in L^1,~i=1,~2.$
We say that the triple of weights $(v,~\omega_1, ~\omega_2)$
satisfies the condition $S_{\overrightarrow{p}},$ if
there exists a positive constant $C$ such that
\be
\big(\int_{\{\tau<\infty\}}\mathcal{M}
(\sigma_1\chi_{\{\tau<\infty\}},\sigma_2\chi_{\{\tau<\infty\}})^{p}vd\mu\big)^{\frac{1}{p}}\leq
C|\{\tau<\infty\}|^{\frac{1}{p_1}}_{\sigma_1}|\{\tau<\infty\}|^{\frac{1}{p_2}}_{\sigma_2},
~\forall \tau\in \mathcal{T}.\ee
\end{definition}

\begin{remark}
If we substitute $p_1=p_2$
and $\omega_1=\omega_2$ into Definition \ref{definition Ap} and Definition \ref{definition Sp},
they reduce to $A_p$ condition and $S_p$ condition in martingale spaces, respectively.\end{remark}
\subsection{Bilinear Version of Two-weight Inequalities}

\begin{theorem} \label{thm_AP}Let $v,~\omega_1, ~\omega_2$ be weights and $1< p_1, ~p_2<\infty.$
Suppose that
$\frac{1}{p}=\frac{1}{p_1 }+\frac{1}{p_2 }$ and $(\omega_1, \omega_2)\in RH(p_1, p_2),$
then the following statements
are equivalent:
\begin{enumerate}
\item  There exists a positive constant $C$ such that
      \begin{equation}\label{Th_B_1}\Big(\int_{\{\tau<\infty\}}|f_{\tau}||g_{\tau}|^pvd\mu\Big)^{\frac{1}{p}}
      \leq C\|f\|_{L^{p_1}(\omega_1)}\|g\|_{L^{p_2}(\omega_2)},
      ~\forall \tau\in\mathcal{T},~ f\in L^{p_1}(\omega_1),~g\in L^{p_2}(\omega_2);\end{equation}
\item  There exists a positive constant $C$ such that
      \begin{equation}\label{Th_B_2}\|\mathcal{M}(f,g)\|_{L^{p,\infty}(v)}\leq
       C\|f\|_{L^{p_1}(\omega_1)}\|g\|_{L^{p_2}(\omega_2)},~
      \forall f\in L^{p_1}(\omega_1),~g\in L^{p_2}(\omega_2);
      \end{equation}
\item The triple of weights $(v,~\omega_1, ~\omega_2)$ satisfies the condition $A_{\overrightarrow{p}},$ i.e.
\begin{equation}\label{Th_B_3}
(v,~\omega_1, ~\omega_2)\in A_{\overrightarrow{p}}.
      \end{equation}
\end{enumerate}
\end{theorem}
\noindent{\bf Proof } We shall follow the scheme:
$(1)\Leftrightarrow(2)\Leftarrow(3)\Leftarrow(1).$

$(1)\Rightarrow(2).$ Let $f\in L^{p_1}(\omega_1),g\in L^{p_2}(\omega_2).$ For
$\lambda>0,$ define $\tau=\inf\{n:|f_n||g_n|>\lambda\}.$ It follows from $\eqref{Th_B_1}$ that
\be\lambda|\{\mathcal{M}(f,g)>\lambda\}|_v^{\frac{1}{p}}
   &=&(\int_{\{\tau<\infty\}}\lambda^p vd\mu )^{\frac{1}{p}}\\
   &\leq&(\int_{\{\tau<\infty\}}|f_{\tau}|^p|g_{\tau}|^pvd\mu )^{\frac{1}{p}}\\
   &\leq&C\|f\|_{L^{p_1}(\omega_1)}\|g\|_{L^{p_2}(\omega_2)}.
   \ee
Thus $\eqref{Th_B_2}$ is valid.

$(2)\Rightarrow(1).$ Fix $n\in N$ and $B\in\mathcal {F}_n.$ For
$f\in L^{p_1}(\omega_1)$ and $g\in L^{p_2}(\omega_2),$ let $$F=f\chi_B \text{ and } G=g\chi_B,$$ respectively.
Then $E_n(F)=f_n\chi_B$ and $E_n(G)=g_n\chi_B.$
Moreover
$$
|f_ng_n|\chi_B\leq\mathcal{M}(F,G).
$$
Combing with $\eqref{Th_B_2},$ we have \be\lambda^p\int_{B\cap\{
        |f_n g_n|>\lambda\}}vd\mu
        &\leq& \lambda^p\int_{\{\mathcal{M}(F,G)>\lambda\}}vd\mu \\
&\leq& C \Big(\int_\Omega|F|^{p_1}\omega_1
        d\mu\Big)^{\frac{p}{p_1}}\Big(\int_\Omega|G|^{p_2}\omega_2
        d\mu\Big)^{\frac{p}{p_2}}\\
&=& C\Big(\int_B|f|^{p_1}\omega_1
        d\mu\Big)^{\frac{p}{p_1}}\Big(\int_B|g|^{p_2}\omega_2
        d\mu\Big)^{\frac{p}{p_2}}.\ee
For $k\in Z$, let
\be
B_k = \{2^k<
    |f_n|| g_n|\leq2^{k+1}\}.
\ee
Note that $$\{2^k<|f_n|| g_n|\leq2^{k+1}\}
\subseteq \{2^k<|f_n||g_n|\},$$
then
\be \int_\Omega
    (|f_n|| g_n|)^pvd\mu
&=& \sum\limits_{k\in Z}\int_{B_k}
    (|f_n|| g_n|)^pvd\mu\\
&\leq&C \sum\limits_{k\in Z}\int_{{B_k}\cap\{
    |f_n|| g_n|>2^k\}}2^{kp}vd\mu\\
&\leq&C \sum\limits_{k\in Z}\Big(\int_{B_k}|f|^{p_1}\omega_1
    d\mu\Big)^{\frac{p}{p_1}}\Big(\int_{B_k}|g|^{p_2}\omega_2
    d\mu\Big)^{\frac{p}{p_2}}\\
&\leq&C \Big(\sum\limits_{k\in Z}\int_{B_k}|f|^{p_1}\omega_1
    d\mu\Big)^{\frac{p}{p_1}}\Big(\sum\limits_{k\in Z}\int_{B_k}|g|^{p_2}\omega_2
    d\mu\Big)^{\frac{p}{p_2}}\\
&=&C \Big(\int_\Omega|f|^{p_1}\omega_1
    d\mu\Big)^{\frac{p}{p_1}}\Big(\int_\Omega|g|^{p_2}\omega_2
    d\mu\Big)^{\frac{p}{p_2}},\ee
where we have used Holder's inequality. As for $\tau\in\mathcal {T},$
it is easy to see that
\be\int_{\{\tau<\infty\}}(|f_{\tau}||g_{\tau}|)^pvd\mu
   &=&\sum\limits_{n\geq0}\int_{\{\tau=n\}}(|f_n||g_n|)^pvd\mu\\
   &\leq&C\sum\limits_{n\geq0}\Big(\int_\Omega|f\chi_{\{\tau=n\}}|^{p_1}\omega_1
         d\mu\Big)^{\frac{p}{p_1}}\Big(\int_\Omega|g\chi_{\{\tau=n\}}|^{p_2}\omega_2
         d\mu\Big)^{\frac{p}{p_2}}\\
   &\leq&C\Big(\sum\limits_{n\geq0}\int_\Omega|f\chi_{\{\tau=n\}}|^{p_1}\omega_1
         d\mu\Big)^{\frac{p}{p_1}}
         \Big(\sum\limits_{n\geq0}\int_\Omega|g\chi_{\{\tau=n\}}|^{p_2}\omega_2
         d\mu\Big)^{\frac{p}{p_2}}\\
   &\leq&C\Big(\int_\Omega|f|^{p_1}\omega_1
         d\mu\Big)^{\frac{p}{p_1}}
         \Big(\int_\Omega|g|^{p_2}\omega_2
         d\mu\Big)^{\frac{p}{p_2}}.\ee
Therefore, $$\Big(\int_{\{\tau<\infty\}}(|f_{\tau}||g_{\tau}|)^pvd\mu\Big)^{\frac{1}{p}}
      \leq C\|f\|_{L^{p_1}(\omega_1)}\|g\|_{L^{p_2}(\omega_2)}.$$

$(3)\Rightarrow(2).$ For
$f\in L^{p_1}(\omega_1),\ g\in L^{p_1}(\omega_2)$ and $n\in N,$ we get
$$|E_n(f)|\leq E_n(|f^{p_1}\omega_1|)^{\frac{1}{p_1}}E_n(\omega_1^{-\frac{1}{p_1-1}})^{\frac{1}{p'_1}}
\text{ and }
|E_n(g)|\leq E_n(|g^{p_2}\omega_2|)^{\frac{1}{p_2}}E_n(\omega_2^{-\frac{1}{p_2-1}})^{\frac{1}{p'_2}}.$$
Furthermore,
\be|E_n(f)E_n(g)|^p
   &\leq&E_n(|f^{p_1}\omega_1|)^{\frac{p}{p_1}}E_n(|g^{p_2}\omega_2|)^{\frac{p}{p_2}}
   E_n(\omega_1^{-\frac{1}{p_1-1}})^{\frac{p}{p'_1}}E_n(\omega_2^{-\frac{1}{p_2-1}})^{\frac{p}{p'_2}}\\
&=&E_n^v(|f^{p_1}\omega_1v^{-1}|)^{\frac{p}{p_1}}E^v_n(|g^{p_2}\omega_2v^{-1}|)^{\frac{p}{p_2}}
   E_n(v)E_n(\omega_1^{-\frac{1}{p_1-1}})^{\frac{p}{p'_1}}
   E_n(\omega_2^{-\frac{1}{p_2-1}})^{\frac{p}{p'_2}},\ee
where $E_n^v(\cdot)$ is the conditional expectation relative to the probability measure $\frac{v}{|\Omega|_v}d\mu.$
Because of $(\ref{Th_B_3}),$ we get $$|E_n(f)E_n(g)|\leq{\bf C} E_n^v(|f^{p_1}\omega_1v^{-1}|)^{\frac{1}{p_1}}
E^v_n(|g^{p_2}\omega_2v^{-1}|)^{\frac{1}{p_2}}.$$
Thus $$\mathcal{M}(f,g)\leq C M^v(f^{p_1}\omega_1v^{-1})^{\frac{1}{p_1}}
M^v(g^{p_2}\omega_2v^{-1})^{\frac{1}{p_2}}.$$
From this, using Holder's inequality for weak spaces, we obtain
\be\|\mathcal{M}(f,g)\|_{L^{p,\infty}(v)}
&\leq&C\|M^v(f^{p_1}\omega_1v^{-1})^{\frac{1}{p_1}}\|_{L^{p_1,\infty}(v)}
      \|M^v(g^{p_2}\omega_2v^{-1})^{\frac{1}{p_2}}\|_{L^{p_2,\infty}(v)}\\
&=&C\|M^v(f^{p_1}\omega_1v^{-1})\|^{\frac{1}{p_1}}_{L^{1,\infty}(v)}
      \|M^v(g^{p_2}\omega_2v^{-1})\|^{\frac{1}{p_2}}_{L^{1,\infty}(v)}\\
&\leq&C\|f^{p_1}\omega_1v^{-1}\|^{\frac{1}{p_1}}_{L^{1}(v)}
       \|g^{p_2}\omega_2v^{-1}\|^{\frac{1}{p_2}}_{L^{1}(v)}\\
&=&C\|f^{p_1}\omega_1\|^{\frac{1}{p_1}}_{L^{1}}
       \|g^{p_2}\omega_2\|^{\frac{1}{p_2}}_{L^{1}}\\
&=&C\|f\|_{L^{p_1}(\omega_1)}
    \|g\|_{L^{p_2}(\omega_2)}.\ee

$(1)\Rightarrow(3).$ For any $n\in N$ and $B\in \mathcal{F}_n,$ set
$f=\omega_1^{-\frac{1}{p_1-1}}\chi_B$ and $g=\omega_2^{-\frac{1}{p_2-1}}\chi_B.$
Then $$\Big(\int_BE_n(\omega_1^{-\frac{1}{p_1-1}})^p
E_n(\omega_2^{-\frac{1}{p_2-1}})^pvd\mu\Big)^{\frac{1}{p}}
      \leq C\Big(\int_\Omega \omega_1^{-\frac{1}{p_1-1}}\chi_Bd\mu\Big)^{\frac{1}{p_1}}
            \Big(\int_\Omega \omega_2^{-\frac{1}{p_2-1}}\chi_Bd\mu\Big)^{\frac{1}{p_2}}.$$
Furthermore,
\begin{equation}\label{Th_B_4}\Big(\int_BE_n(\omega_1^{-\frac{1}{p_1-1}})^p
E_n(\omega_2^{-\frac{1}{p_2-1}})^pE_n(v)d\mu\Big)^{\frac{1}{p}}
      \leq C\Big(\int_BE_n(\omega_1^{-\frac{1}{p_1-1}})d\mu\Big)^{\frac{1}{p_1}}
            \Big(\int_BE_n(\omega_2^{-\frac{1}{p_2-1}})d\mu\Big)^{\frac{1}{p_2}}.\end{equation}

Thus, there exists a constant $C$ such that
$$\Big(E_n(\omega_1^{-\frac{1}{p_1-1}})^p
E_n(\omega_2^{-\frac{1}{p_2-1}})^pE_n(v)\Big)^{\frac{1}{p}}
      \leq CE_n(\omega_1^{-\frac{1}{p_1-1}})^{\frac{1}{p_1}}
            E_n(\omega_2^{-\frac{1}{p_2-1}})^{\frac{1}{p_2}}.$$
Otherwise, for any $C>0,$ let $$B=\{E_n(\omega_1^{-\frac{1}{p_1-1}})^p
E_n(\omega_2^{-\frac{1}{p_2-1}})^pE_n(v)
      > CE_n(\omega_1^{-\frac{1}{p_1-1}})^{\frac{p}{p_1}}
            E_n(\omega_2^{-\frac{1}{p_2-1}})^{\frac{p}{p_2}}\},$$
then $\mu(B)>0.$ Consequently,
\begin{eqnarray}
\int_BE_n(\omega_1^{-\frac{1}{p_1-1}})^p
        E_n(\omega_2^{-\frac{1}{p_2-1}})^pE_n(v)d\mu
      &>& C\int_BE_n(\omega_1^{-\frac{1}{p_1-1}})^{\frac{p}{p_1}}
             E_n(\omega_2^{-\frac{1}{p_2-1}})^{\frac{p}{p_2}}d\mu\nonumber\\
      &\geq& C\int_BE_n(\omega_1^{-\frac{1}{p_1-1}\frac{p}{p_1}}
             \omega_2^{-\frac{1}{p_2-1}\frac{p}{p_2}})d\mu\label{Holder}\\
      &=& C\int_B\omega_1^{-\frac{1}{p_1-1}\frac{p}{p_1}}
             \omega_2^{-\frac{1}{p_2-1}\frac{p}{p_2}}d\mu\nonumber\\
      &\geq& C\big(\int_B\omega_1^{-\frac{1}{p_1-1}}d\mu\big)^{\frac{p}{p_1}}
            \big(\int_B\omega_2^{-\frac{1}{p_2-1}}d\mu\big)^{\frac{p}{p_2}}\label{RHolder},\end{eqnarray}
where $(\ref{Holder})$ and $(\ref{RHolder})$ use Holder's inequality for $E_n(\cdot)$ and $RH(p_1,p_2)$ condition,
respectively. It follows that
\be\int_BE_n(\omega_1^{-\frac{1}{p_1-1}})^pE_n(\omega_2^{-\frac{1}{p_2-1}})^pE_n(v)d\mu> C\big(\int_B\omega_1^{-\frac{1}{p_1-1}}d\mu\big)^{\frac{p}{p_1}}
            \big(\int_B\omega_2^{-\frac{1}{p_2-1}}d\mu\big)^{\frac{p}{p_2}},\ee
which contradicts $(\ref{Th_B_4})$. By contradiction, we have
$$\Big(E_n(\omega_1^{-\frac{1}{p_1-1}})^p
E_n(\omega_2^{-\frac{1}{p_2-1}})^pE_n(v)\Big)^{\frac{1}{p}}
      \leq CE_n(\omega_1^{-\frac{1}{p_1-1}})^{\frac{1}{p_1}}
            E_n(\omega_2^{-\frac{1}{p_2-1}})^{\frac{1}{p_2}}.$$
Then
\be E_n(v)^{\frac{1}{p}}E_n(\omega_1^{1-p'_1})^{\frac{1}{p'_1}}
E_n(\omega_2^{1-p'_2})^{\frac{1}{p'_2}}\leq C.\ee

\begin{theorem}\label{thm_Sp} Let $v,\omega_1, \omega_2$ be weights and $1<p_1, p_2<\infty.$
Suppose that
$\frac{1}{p}=\frac{1}{p_1 }+\frac{1}{p_2 }$ and $(\omega_1, \omega_2)\in RH(p_1, p_2),$
then the following statements
are equivalent:\begin{enumerate}
\item There exists a positive constant $C$ such that
\be
\|\mathcal{M}(f,g)\|_{L^p(v)}\leq
C\|f\|_{L^{p_1}(\omega_1)}\|g\|_{L^{p_2}(\omega_2)},
~\forall f\in L^{p_1}(\omega_1),~g\in L^{p_2}(\omega_2);
\ee
\item There exists a positive constant $C$ such that
\begin{equation}
\label{Th_A_1}\|\mathcal{M}(f\sigma_1,g\sigma_2)\|_{L^p(v)}\leq
C\|f\|_{L^{p_1}(\sigma_1)}\|g\|_{L^{p_2}(\sigma_2)},
~\forall f\in L^{p_1}(\sigma_1),~g\in L^{p_2}(\sigma_2),
\end{equation}
where $\sigma_i=\omega_i^{\frac{1}{p_i-1}}, ~i=1,~2;$
\item The triple of weights $(v,~\omega_1, ~\omega_2)$ satisfies the condition $S_{\overrightarrow{p}},$ i.e.
\be
(v,~\omega_1, ~\omega_2)\in S_{\overrightarrow{p}}.
\ee  \end{enumerate}
\end{theorem}

\begin{remark}We mention that the first author has also
obtained a similar characterization for the multisublinear maximal
function in function space.
The multilinear testing condition was
further discussed by \cite{Chen-Damian} in function space,
which generalized the result in \cite{Sawyer E T.}.
\end{remark}

\noindent{\bf Proof } It is clear that $(1)\Leftrightarrow(2)\Rightarrow(3),$
so we omit them.  To prove $(3)\Rightarrow(2),$ we proceed in the following way.
Let $f\in L^{p_1}(\sigma_1),g\in L^{p_2}(\sigma_2)$. For all $k\in Z$, define
stopping times $$\tau_k=\inf\{n:|E(f\sigma_1|\mathcal{F}_n)E(g\sigma_2|
\mathcal{F}_n)|>2^k\}.$$
Set
\be
A_{k,j}=\{\tau_k<\infty\}\cap\{2^j<E(\sigma_1|\mathcal{F}_{\tau_k})E(\sigma_2|
\mathcal{F}_{\tau_k})\leq2^{j+1}\};
\ee
\be
B_{k,j}=\{\tau_k<\infty,\tau_{k+1}=\infty\}\cap\{2^j<E(\sigma_1|\mathcal{F}_{\tau_k})
E(\sigma_2|\mathcal{F}_{\tau_k})\leq2^{j+1}\},~ j\in Z.
\ee
Then $A_{k,j}\in \mathcal{F}_{\tau_k}, B_{k,j}\subseteq A_{k,j}$.
Moreover, $\{B_{k,j}\}_{k,j}$ is a family of disjoint sets and
\be
\{2^k<\mathcal{M}(f\sigma_1,g\sigma_2)\leq2^{k+1}\}=\{\tau_k<\infty,\tau_{k+1}=\infty\}=\bigcup\limits_{j\in
Z} B_{k,j}, k\in Z.
\ee
Trivially,
\be
E(f\sigma_1|\mathcal{F}_{\tau_k})=E^{\sigma_1}(f|\mathcal{F}_{\tau_k})E(\sigma_1|\mathcal{F}_{\tau_k})\text{ and }
E(g\sigma_2|\mathcal{F}_{\tau_k})=E^{\sigma_2}(g|\mathcal{F}_{\tau_k})E(\sigma_2|\mathcal{F}_{\tau_k}).
\ee
On each $A_{k,j},$ we have \be 2^{kp}&\leq&
                 \mathop{\hbox{ess inf}}\limits_{A_{k,j}}
                 |E(f\sigma_1|\mathcal{F}_{\tau_k})^pE(g\sigma_2|\mathcal{F}_{\tau_k})^p|\\
              &\leq&\mathop{\hbox{ess inf}}\limits_{A_{k,j}}
                 |E^{\sigma_1}(f|\mathcal{F}_{\tau_k})E^{\sigma_2}(g|\mathcal{F}_{\tau_k})|^p
                 \mathop{\hbox{ess sup}}\limits_{A_{k,j}}
                 \big(E(\sigma_1|\mathcal{F}_{\tau_k})E(\sigma_2|\mathcal{F}_{\tau_k})\big)^p\\
              &\leq&2^p\mathop{\hbox{ess inf}}\limits_{A_{k,j}}
                 |E^{\sigma_1}(f|\mathcal{F}_{\tau_k})E^{\sigma_2}(g|\mathcal{F}_{\tau_k})|^p
                 |B_{k,j}|_v^{-1}\int_{B_{k,j}}
                 \big(E(\sigma_1|\mathcal{F}_{\tau_k})E(\sigma_2|\mathcal{F}_{\tau_k})\big)^pvd\mu.\ee
To estimate $\int_\Omega\mathcal{M}(f\sigma_1,g\sigma_2)^p v d\mu, $ firstly we have \be
   &~&\int_\Omega\mathcal{M}(f\sigma_1,g\sigma_2)^p v d\mu\\
   &=&\sum\limits_{k\in Z}\int_{\{2^k<\mathcal{M}(f\sigma_1,g\sigma_2)\leq2^{k+1}\}}
      \mathcal{M}(f\sigma_1,g\sigma_2)^p v d\mu\\
   &\leq&2^p\sum\limits_{k\in Z}\int_{\{2^k<\mathcal{M}(f\sigma_1,g\sigma_2)\leq2^{k+1}\}}2^{kp}vd\mu\\
   &=&2^p\sum\limits_{k\in Z,j\in Z}2^{kp}\int_{B_{k,j}} v d\mu\\
   &\leq&4^p\sum\limits_{k\in Z,j\in Z}\mathop{\hbox{ess inf}}\limits_{A_{k,j}}
      |E^{\sigma_1}(f|\mathcal{F}_{\tau_k})E^{\sigma_2}(g|\mathcal{F}_{\tau_k})|^p
      \int_{B_{k,j}}
      \big(E(\sigma_1|\mathcal{F}_{\tau_k})E(\sigma_2|\mathcal{F}_{\tau_k})\big)^pvd\mu.
\ee
It is clear that $\vartheta$ is a measure on $X=Z^2$ with
\be
\vartheta(k,j)=\int_{B_{k,j}}
      \big(E(\sigma_1|\mathcal{F}_{\tau_k})E(\sigma_2|\mathcal{F}_{\tau_k})\big)^pvd\mu.
\ee
For the above $f\in L^{p_1}(\sigma_1),\ g\in L^{p_2}(\sigma_2),$ define \be
T_{f,g}(k,j)=\mathop{\hbox{ess inf}}\limits_{A_{k,j}}
      |E^{\sigma_1}(f|\mathcal{F}_{\tau_k})E^{\sigma_2}(g|\mathcal{F}_{\tau_k})|^p\ee
and denote
\be E_\lambda=\Big\{(k,j):
\mathop{\hbox{ess inf}}\limits_{A_{k,j}}
|E^{\sigma_1}(f|\mathcal{F}_{\tau_k})E^{\sigma_2}(g|\mathcal{F}_{\tau_k})|^p>\lambda\Big\}
\hbox{ and }G_\lambda=\bigcup\limits_{(k,j)\in E_\lambda}A_{k,j} \ee
for each $\lambda>0. $ Then we have \be |\{T_{f,g}>\lambda\}|_\vartheta
&=&\sum\limits_{(k,j)
        \in E_\lambda}\int_{B_{k,j}}
        \big(E(\sigma_1|\mathcal{F}_{\tau_k})E(\sigma_2|\mathcal{F}_{\tau_k})\big)^pvd\mu\\
&\leq&\sum\limits_{(k,j)
        \in E_\lambda}\int_{B_{k,j}}\big(E(\sigma_1\chi_{G_\lambda}|\mathcal{F}_{\tau_k})
        E(\sigma_2\chi_{G_\lambda}|\mathcal{F}_{\tau_k})\big)^pvd\mu\\
&\leq& \int_{G_\lambda}\mathcal{M}(\sigma_1\chi_{G_\lambda},\sigma_2\chi_{G_\lambda})^p v d\mu. \ee
Let $\tau=\inf\Big\{n:~
          |E^{\sigma_1}(f|\mathcal{F}_n)E^{\sigma_2}
          (g|\mathcal{F}_n)|^p>\lambda\Big\},$ we have
$G_\lambda\subseteq\Big\{\mathcal{M}^{\sigma_1,\sigma_2}(f,g)^p>\lambda\Big\}=\{\tau<\infty\}.$
It follows from $S_{\overrightarrow{p}}$ and $RH(p_1,p_2)$ that
 \be |\{T_{f,g}>\lambda\}|_\vartheta
&\leq&\int_{\{\tau<\infty\}}\mathcal{M}(\sigma_1\chi_{\{\tau<\infty\}},\sigma_2\chi_{\{\tau<\infty\}})^p v d\mu.\\
&\leq& C|\{\tau<\infty\}|^{\frac{p}{p_1}}_{\sigma_1}|\{\tau<\infty\}|^{\frac{p}{p_2}}_{\sigma_2}\\
&\leq& C\int_{\{\tau<\infty\}}\sigma_1^{\frac{p}{p_1}}\sigma_2^{\frac{p}{p_2}}d\mu.
 \ee
Therefore,
\be~~~~~~\int_\Omega\mathcal{M}(f\sigma_1,g\sigma_2)^p v d\mu
&\leq&4^p\int_XT_{f,g}d\vartheta=4^p\int_0^\infty|\{T_{f,g}>\lambda\}|_\vartheta
              d\lambda\\
&\leq&C\int_0^\infty\int_{\{\tau<\infty\}}\sigma_1^{\frac{p}{p_1}}
              \sigma_2^{\frac{p}{p_2}}d\mu d\lambda \\
&=&C\int_0^\infty\int_{\{\mathcal{M}^{\sigma_1,\sigma_2}(f,g)^p>\lambda\}}\sigma_1^{\frac{p}{p_1}}
              \sigma_2^{\frac{p}{p_2}}d\mu d\lambda \\
&=&C\int_\Omega\mathcal{M}^{\sigma_1,\sigma_2}(f,g)^p\sigma_1^{\frac{p}{p_1}}
              \sigma_2^{\frac{p}{p_2}}d\mu\\
&\leq&C\int_\Omega M^{\sigma_1}(f)^pM^{\sigma_2}(g)^p\sigma_1^{\frac{p}{p_1}}
              \sigma_2^{\frac{p}{p_2}}d\mu\\
&\leq&C\big(\int_\Omega M^{\sigma_1}(f)^{p_1}\sigma_1d\mu\big)^{\frac{p}{p_1}}
              \big(\int_\Omega M^{\sigma_1}(f)^{p_2}\sigma_2d\mu\big)^{\frac{p}{p_2}}\\
&\leq&C\|f\|_{L^{p_1}(\sigma_1)}^{p}\|g\|_{L^{p_2}(\sigma_2)}^{p}.
\ee where we have used Holder's inequality.
Whence $(\ref{Th_A_1})$ is valid.

\begin{cor} \label{cor_AP1}Let $v,~\omega$ be weights and $1< p<\infty.$
Suppose that $\omega^{-\frac{1}{p-1}}\in L^1.$ Then the following statements
are equivalent:
\begin{enumerate}
\item  There exists a positive constant $C$ such that
      \be\Big(\int_{\{\tau<\infty\}}|f_{\tau}|^pvd\mu\Big)^{\frac{1}{p}}
      \leq C\|f\|_{L^{p}(\omega)},
      ~\forall \tau\in\mathcal{T},~ f\in L^{p}(\omega);\ee
\item  There exists a positive constant $C$ such that
      \be\|M f\|_{L^{p,\infty}(v)}\leq
       C\|f\|_{L^{p}(\omega)},~
      \forall f\in L^{p}(\omega);
      \ee
\item The couple of weights $(v,~\omega)$ satisfies the condition $A_{p},$ i.e.
\be
(v,~\omega)\in A_p.
      \ee
\end{enumerate}
\end{cor}

\begin{cor}\label{cor_Sp1} Let $v,\omega$ be weights and $1<p<\infty.$ Suppose that $\omega^{-\frac{1}{p-1}}\in L^1.$
Then the following statements
are equivalent:\begin{enumerate}
\item There exists a positive constant $C$ such that
\be
\|Mf\|_{L^p(v)}\leq
C\|f\|_{L^{p}(\omega},
~\forall f\in L^{p}(\omega);
\ee
\item There exists a positive constant $C$ such that
\be\|M(f\sigma)\|_{L^p(v)}\leq
C\|f\|_{L^{p}(\sigma)},
~\forall f\in L^{p}(\sigma),
\ee
where $\sigma=\omega^{\frac{1}{p-1}};$
\item The couple of weights $(v,~\omega)$ satisfies the condition $S_p,$ i.e.
\be
(v,~\omega)\in S_p.
\ee  \end{enumerate}
\end{cor}

\noindent{\bf Proof } If we substitute $p_1=p_2$
and $\omega_1=\omega_2$ into Theorem \ref{thm_AP} and Theorem \ref{thm_Sp},
then condition $RH$ is trivial and we get
Corollary \ref{cor_AP1} and Corollary \ref{cor_Sp1}.

\subsection{Bilinear version of One-weight Theory}
Firstly, we recall the following Proposition \ref{prop Ap} for $A_p$ weight in martingale context(\cite{M. Kikuchi, R. L. Long}).
Then, we partially give its bilinear analogue.

\begin{prop}
\label{prop Ap}Let $\omega$ be weights and $1< p<\infty.$
Suppose that $\omega^{-\frac{1}{p-1}}\in L^1.$ Then the following statements
are equivalent:
\begin{enumerate}
\item The weight $\omega$ satisfies the condition $A_{p},$ i.e.
      \be
\sup\limits_{n\geq0}E_n(\omega)E_n(\omega^{-\frac{1}{p-1}})^{p-1}
      \leq C;
      \ee
\item  There exists a positive constant $C$ such that
      $\|E_n\|_{L^p(\omega)\rightarrow L^p(\omega)}
      \leq C~, \forall n \in N,$ i.e.
      \be\|E_n(f)\|_{L^{p}(\omega)}
      \leq C\|f\|_{L^{p}(\omega)},
      ~\forall n\in N,~ f\in L^{p}(\omega);\ee
\item  There exists a positive constant $C$ such that
      \be\lim\limits_{n\rightarrow\infty}
      \Big(\int_\Omega |E_n(f)-f|^p\omega d\mu\Big)^{\frac{1}{p}}=0,
          ~\forall f\in L^{p}(\omega);
      \ee
\item There exists a positive constant $C$ such that
      \be
           \|Mf\|_{L^p(\omega)}\leq
          C\|f\|_{L^{p}(\omega)},
          ~\forall f\in L^{p}(\omega).
      \ee
\end{enumerate}
\end{prop}

\begin{prop} \label{cor_AP2}Let $v,~\omega_1, ~\omega_2$ be weights and $1< p_1, ~p_2<\infty.$
Suppose that
$\frac{1}{p}=\frac{1}{p_1 }+\frac{1}{p_2 }$ and $(\omega_1, \omega_2)\in RH(p_1, p_2),$
then the following statements
are equivalent:
\begin{enumerate}
\item  There exists a positive constant $C$ such that
      \begin{equation}\label{cor_B_11}\Big(\int_\Omega |E_n(f)E_n(g)|^pvd\mu\Big)^{\frac{1}{p}}
      \leq C\|f\|_{L^{p_1}(\omega_1)}\|g\|_{L^{p_2}(\omega_2)},
      ~\forall n\in N,~ f\in L^{p_1}(\omega_1),~g\in L^{p_2}(\omega_2);\end{equation}
\item The triple of weights $(v,~\omega_1, ~\omega_2)$ satisfies
the condition $A_{\overrightarrow{p}},$ i.e.
\be
(v,~\omega_1, ~\omega_2)\in A_{\overrightarrow{p}}.
      \ee
\end{enumerate}
\end{prop}

\noindent{\bf Proof } It is easy to check that (\ref{cor_B_11})
and (\ref{Th_B_1}) are equivalent. Thus we get
Proposition \ref{cor_AP2} from Theorem \ref{thm_AP}.

\begin{remark}
The condition $(\omega_1, \omega_2)\in RH(p_1, p_2)$ is only used to prove that (\ref{cor_B_11}) implies
the condition $(v,~\omega_1, ~\omega_2)\in A_{\overrightarrow{p}}$ in Proposition \ref{cor_AP2}.
\end{remark}

\begin{lemma} \label{lem con}Let $\omega_1, ~\omega_2$ be weights and $1<p_1, ~p_2<\infty.$
Suppose that $\frac{1}{p}=\frac{1}{p_1 }+\frac{1}{p_2},$
$\omega_i^{-\frac{1}{p_i-1}}\in L^1,$ $i=1,~2$ and $v=\omega_1^{\frac{p}{p_1}}\omega_2^{\frac{p}{p_2}}.$
If $f\in L^{p_1}(\omega_1),~g\in L^{p_2}(\omega_2)$ and
$E_n(f)E_n(g)\in L^p(v),~\forall n\in N,$
then \begin{equation}\label{lem cov}\lim\limits_{n\rightarrow\infty}
      \Big(\int_\Omega |E_n(f)E_n(g)-fg|^pvd\mu\Big)^{\frac{1}{p}}=0,\end{equation}
if and only if, for any $\varepsilon>0,$ there is a nonnegative function $y \in L^p(v)$ such that
\begin{equation}\label{lem cov1}\sup\limits_{n\geq0}
      \Big(\int_\Omega |E_n(f)E_n(g)\chi_{\{|E_n(f)E_n(g)|\geq y\}}|^pvd\mu\Big)^{\frac{1}{p}}
      \leq\varepsilon.\end{equation}
\end{lemma}

\noindent{\bf Proof } Suppose that (\ref{lem cov1}) is valid, we will
prove (\ref{lem cov}). For any $\varepsilon>0,$ there is a nonnegative function $y \in L^p(v)$ such that
$$\sup\limits_{n\geq0}
      \Big(\int_\Omega |E_n(f)E_n(g)\chi_{\{|E_n(f)E_n(g)|\geq y\}}|^pvd\mu\Big)^{\frac{1}{p}}
      \leq\varepsilon.$$ Since $\|fg\|_{L^p(v)}\leq \|f\|_{L^{p_1}(\omega_1)}\|g\|_{L^{p_2}(\omega_2)}<\infty,$
we can assume that $y>|fg|.$
We also have $\lim \limits_{n\rightarrow \infty}f_n=f$ and $\lim \limits_{n\rightarrow \infty}g_n=g,$
because martingale $(f_n)_{n\geq0}$ and martingale $(g_n)_{n\geq0}$ are uniformly integral. Thus
$$(2y)^p\geq |f_ng_n\chi_{\{|f_ng_n|<y\}}-fg|^p\rightarrow0,~\hbox{as}~n\rightarrow \infty.$$
It follows from dominated integral theorem that
$$\lim \limits_{n\rightarrow \infty}\|f_ng_n\chi_{\{|f_ng_n|<y\}}-fg\|_{L^p(v)}=0.$$
For the above $\varepsilon,$
there is a $n_0\in N,$ such that
$$\|f_ng_n\chi_{\{|f_ng_n|<y\}}-fg\|_{L^p(v)}<\varepsilon,~\forall n>n_0.$$ Moreover,
\be\|f_ng_n-fg\|_{L^p(v)}
     &=&\|f_ng_n(\chi_{\{|f_ng_n|<y\}}+
           \chi_{\{|f_ng_n|\geq y\}})-fg\|_{L^p(v)}\\
     &\leq&(2^{\frac{1-p}{p}}\vee1)\big(\|f_ng_n\chi_{\{|f_ng_n|<y\}}-fg\|_{L^p(v)}+
           \|f_ng_n\chi_{\{|f_ng_n|\geq y\}}\|_{L^p(v)}\big)\\
     &<&2(2^{\frac{1-p}{p}}\vee1)\varepsilon,~\forall n>n_0,\ee
which implies (\ref{lem cov}).

     Conversely, we assume that (\ref{lem cov}) is valid. Since $fg\in L^p(v),$ we obtain that
for any $0<\varepsilon<1,$ there exists $\delta>0$ such that
whenever $E\in\mathcal{F}$ satisfies $|E|_v<\delta,$ then
$\Big(\int_E |fg|^pvd\mu\Big)^{\frac{1}{p}}<\frac{1}{2(2^{\frac{1-p}{p}}\vee1)}\varepsilon.$
For the above $\varepsilon>0,$ there exists $n_0,$ such that
$$\Big(\int_\Omega |E_n(f)E_n(g)-fg|^pvd\mu\Big)^{\frac{1}{p}}<
\big(\frac{1}{2(2^{\frac{1-p}{p}}\vee1)}\wedge\delta^{\frac{1}{p}}\big)\varepsilon,~\forall n\geq n_0.$$
Moreover, for the above $\varepsilon>0,~n\geq n_0,$ we obtain that
\be |\{|E_n(f)E_n(g)|-|fg|>\varepsilon\}|_v
&=&\frac{1}{\varepsilon^p}\int_{\{|E_n(f)E_n(g)|-|fg|>\varepsilon\}}\varepsilon^pvd\mu\\
&\leq&\frac{1}{\varepsilon^p}\int_\Omega |E_n(f)E_n(g)-fg|^pvd\mu<\delta.\ee
Let $y=\max\{2|f_1g_1|,~2|f_2g_2|,~\cdot\cdot\cdot,2|f_{n_0}g_{n_0}|,~|fg|+2\varepsilon\},$
it follows that  $y\in L^p(v)$ and
\be&~&\sup\limits_{n\geq 0}
      \Big(\int_\Omega |E_n(f)E_n(g)\chi_{\{|E_n(f)E_n(g)|\geq y\}}|^pvd\mu\Big)^{\frac{1}{p}}\\
&=&\sup\limits_{n>n_0}
      \Big(\int_{\{|E_n(f)E_n(g)|\geq y\}} |E_n(f)E_n(g)|^pvd\mu\Big)^{\frac{1}{p}}\\
&=&\sup\limits_{n>n_0}
      \Big(\int_{\{|E_n(f)E_n(g)|\geq y\}} |E_n(f)E_n(g)-fg
      +fg|^pvd\mu\Big)^{\frac{1}{p}}\\
&\leq&(2^{\frac{1-p}{p}}\vee1)\sup\limits_{n>n_0}
      \Big(\int_\Omega |E_n(f)E_n(g)-fg|^pvd\mu\Big)^{\frac{1}{p}}+
      \\&~&+
      (2^{\frac{1-p}{p}}\vee1)\sup\limits_{n>n_0}
      \Big(\int_{\{|E_n(f)E_n(g)|-|fg|>\varepsilon\}}|f g|^pvd\mu\Big)^{\frac{1}{p}}\\
&<&\varepsilon. \ee
We are done.

\begin{prop}\label{prop_AP1}Let $\omega_1, ~\omega_2$ be weights and $1< p_1, ~p_2<\infty.$
Suppose that
$\frac{1}{p}=\frac{1}{p_1 }+\frac{1}{p_2 }$
and $v=\omega_1^{\frac{p}{p_1}}\omega_2^{\frac{p}{p_2}}.$
If the triple of weights $(v,~\omega_1, ~\omega_2)$ satisfies
the condition $A_{\overrightarrow{p}},$ then \begin{equation}\label{prop cov}\lim\limits_{n\rightarrow\infty}
      \Big(\int_\Omega |E_n(f)E_n(g)-fg|^pvd\mu\Big)^{\frac{1}{p}}=0,
~\forall f\in L^{p_1}(\omega_1),~g\in L^{p_2}(\omega_2).\end{equation}
\end{prop}

\noindent{\bf Proof } Let $f\in L^{p_1}(\omega_1),~g\in L^{p_2}(\omega_2).$
It follows from condition $A_{\overrightarrow{p}}$ and Proposition \ref{cor_AP2}
that $$\Big(\int_\Omega |E_n(f)E_n(g)|^pvd\mu\Big)^{\frac{1}{p}}
      \leq C\|f\|_{L^{p_1}(\omega_1)}\|g\|_{L^{p_2}(\omega_2)},
      ~\forall n\in N,$$
which is the assumption of the Lemma \ref{lem con}. If
(\ref{lem cov1}) is valid , we have (\ref{prop cov}) by the Lemma \ref{lem con}. We
will prove (\ref{lem cov1}) in the following way. Since $f$ and $g$ are integral,
 martingale $(f_n)_{n\geq0}$ and martingale $(g_n)_{n\geq0}$ are uniformly integral.
It follows from Doob's inequality that
\begin{equation}\label{prop_ref2}\sup\limits_{\lambda>0}\lambda|\{Mf>\lambda\}|\leq \int_{\Omega}|f|d\mu
\hbox{ and }\sup\limits_{\lambda>0}\lambda|\{Mg>\lambda\}|\leq \int_{\Omega}|g|d\mu.\end{equation}
For $n\in N,$ fix $\lambda>0$ which will be determined later. Then,
\begin{eqnarray}
&~&\Big(\int_\Omega |E_n(f)E_n(g)
          \chi_{\{|E_n(f)E_n(g)|\geq \lambda\}}|^pvd\mu\Big)^{\frac{1}{p}}\nonumber\\
&=&\Big(\int_\Omega |E_n(f\chi_{\{|E_n(f)E_n(g)|\geq \lambda\}})
          E_n(g\chi_{\{|E_n(f)E_n(g)|\geq \lambda\}})|^pvd\mu\Big)^{\frac{1}{p}}\nonumber\\
&\leq&\Big(\int_\Omega E_n(|f\chi_{\{MfMg\geq\lambda\}}|)
          E_n(|g\chi_{\{MfMg\geq \lambda\}}|)^pvd\mu\Big)^{\frac{1}{p}}\nonumber\\
&\leq&\label{prop_ref}C\|f\chi_{\{MfMg\geq \lambda\}}\|_{L^{p_1}(\omega_1)}
          \|g\chi_{\{MfMg\geq \lambda\}}\|_{L^{p_2}(\omega_2)},\end{eqnarray}
where (\ref{prop_ref}) is a result of Corollary \ref{cor_AP2}. It is clear that
$$\{MfMg\geq \lambda\}\subseteq\{Mf\geq \lambda^{\frac{p}{p_1}}\}\cup\{Mg\geq \lambda^{\frac{p}{p_2}}\}.$$
Thus $|\{MfMg\geq \lambda\}|\leq|\{Mf\geq \lambda^{\frac{p}{p_1}}\}|+|\{Mg\geq \lambda^{\frac{p}{p_2}}\}|.$
Combing with (\ref{prop_ref2}), we get $\lim\limits_{\lambda\rightarrow\infty}|\{MfMg\geq \lambda\}|=0.$
Then, (\ref{lem cov1}) follows from (\ref{prop_ref}), because of the absolute continuity of the integral.

\begin{prop}\label{prop_bou-cov}Let $\omega_1, ~\omega_2$ be weights and $1<p_1, ~p_2<\infty.$
Suppose that
$\frac{1}{p}=\frac{1}{p_1 }+\frac{1}{p_2 }$
and $v=\omega_1^{\frac{p}{p_1}}\omega_2^{\frac{p}{p_2}}.$
If there exists a positive constant $C$ such that
\be
\|\mathcal{M}(f,g)\|_{L^p(v)}\leq
C\|f\|_{L^{p_1}(\omega_1)}\|g\|_{L^{p_2}(\omega_2)},
~\forall f\in L^{p_1}(\omega_1),~g\in L^{p_2}(\omega_2),\ee
we have $(v,~\omega_1, ~\omega_2)\in A_{\overrightarrow{p}},$ (\ref{cor_B_11}) as well as (\ref{prop cov}).
\end{prop}

\begin{remark}
The proof of Proposition \ref{prop_bou-cov} is clear and we omit it. But we can not give the reverse of the Proposition \ref{prop_bou-cov} in martingale spaces.
\end{remark}
%
%


\begin{thebibliography}{999}

\bibitem {X. Q. Chang} X. Q. Chang,
\newblock {\em {Some Sawyer type inequalities for martingales}},
\newblock Studia Math. {\bf{111}}(1994), 187-194.

\bibitem {Chen-Damian}
W. Chen and W. Damian,
\newblock {\em {Weighted estimates for the multisublinear maximal function}},
\newblock Rend. Circ. Mat. Palermo (2) 62(2013), 379-391.

\bibitem {W. Chen}
W. Chen and P. D. Liu,
\newblock {\em {Weighted inequalities for the generalized maximal operator in martingale spaces}},
\newblock Chin. Ann. Math. {\bf{32}}(2011), 781-792.

\bibitem{Cruz-Uribe D.} D. Cruz-Uribe,
\newblock {\em {The minimal operator and  the geometric maximal operator in $R^n$}},
\newblock Studia Math. {\bf{144}}(2001), 1-37.

\bibitem{Cruz-Uribe D. J. M. Martell}
D. Cruz-Uribe, J. M. Martell and C. Perez,
\newblock {\em{Weights, extrapolation and the theory of Rubio de Francia}},
\newblock Birkhauser/Springer Basel AG, 2011, Basel.

\bibitem{Garcia Rubio}
J. Garcia-Cuerva and J. L. Rubio de Francia,
\newblock {\em {Weighted Norm Inequalities and Related Topics}},
\newblock Amsterdam: North Holland, 1985.

\bibitem{Grafakos L. Liu L. G. Perez C. Torres R. H.}
L. Grafakos, L. G. Liu, C. Perez and R. H. Torres,
\newblock {\em The multilinear strong maximal function},
\newblock J. Geom. Anal. {\bf {21}}(2011), 118-149.

\bibitem{Grafakos L. Liu L. G. Yang D. C.}
L. Grafakos, L. G. Liu and D. C. Yang,
\newblock {\em Multiple-weighted norm inequalities for maximal multi-linear singular integrals with non-smooth kernels},
\newblock Proc. Roy. Soc. Edinburgh Sect. A. {\bf {141 }}(2011), 755-775.

\bibitem{Hruscev S. V.} S. V. Hruscev,
\newblock  {\em{A description of weights satisfying the $A_\infty$ condition of Muckenhoupt}},
\newblock Proc. Amer. Math. Soc. {\bf{90}}(1984), 253-257.

\bibitem{Wei H. Shi X. L. and Sun Q. Y.} W. Hu, X. L. Shi and Q. Y. Sun,
\newblock {\em{$A_\infty$ condition characterized by maximal geometric mean operator}},
\newblock Lecture Notes in Math. vol. {\bf{1494}}, Springer-Verlag, New York,
1991, 68-72.

\bibitem {M. Izumisawa} M. Izumisawa and N. Kazamaki,
\newblock {\em {Weighted norm inequalities for martingale}},
\newblock Tohoku Math. J. {\bf{29}}(1977), 115-124.

\bibitem {M. Kikuchi} M. Kikuchi,
\newblock {\em {A note on the convergence of martingales in Banach function spaces}},
\newblock Anal. Math. {\bf{25}}(1999), 265-276.

\bibitem{Kurtz D. S.} D. S. Kurtz,
\newblock  {\em Classical operators on mixed-normed spaces with
product weights},
\newblock Rocky Mountain J. Math. {\bf {37}}(2007), 269-283.

\bibitem{Lerner A.K. Ombrosi S.} A. K. Lerner, S. Ombrosi, C. Perez, R. H. Torres and R. Trujillo-Gonzalez,
\newblock  {\em New maximal functions and
multiple weights for the multilinear Calderon-Zygmund theory},
\newblock Adv. Math. {\bf {220}}(2009), 1222-1264.

\bibitem {R. L. Long} R. L. Long,
\newblock {\em {Martingale spaces and inequalities}},
\newblock Peking University Press, 1993, Beijing.

\bibitem{Long R. L. and Peng L. Z.} R. L. Long and L. Z. Peng,
\newblock {\em {Two weighted maximal $(p,q)$ inequalities in martingale setting}},
\newblock Acta Math. Sinica. {\bf{29}}(1986), 253-258.

\bibitem{Muckenhoupt B} B. Muckenhoupt,
\newblock  {\em{Weighted norm inequalities
for the Hardy maximal function}},
\newblock Trans. Amer. Math. Soc. {\bf{165}}(1972), 207-226.

\bibitem{Nehari} Z. Nehari,
\newblock  {\em{Inverse Holder inequalities}},
\newblock J. Math. Anal. Appl. {\bf{21}}(1968), 405-420.

\bibitem {F. Ruiza} F. J. Ruiz,
\newblock {\em {A unified approach to Carleson measures and $A_p$ weights}},
\newblock Pacific J. Math. {\bf{117}}(1985), 397-404.

\bibitem {F. Ruiz} F. J. Ruiz and J. L. Torrea,
\newblock {\em {A unified approach to Carleson measures and $A_p$ weights II}},
\newblock Pacific J. Math. {\bf{120}}(1985), 189-197.

\bibitem {Sawyer E T.} E. T. Sawyer,
\newblock {\em {A characterization of a two weight norm inequality for maximal operators}},
\newblock Studia Math. {\bf{75}}(1982), 1-11.

\bibitem{Sbordone C. and Wik I.} C. Sbordone and I. Wik,
\newblock {\em {Maximal functions and related weight classes}},
\newblock Publ. Mat. {\bf{38}}(1994), 127-155.

\bibitem{Yin X. Q. and Muckenhoupt B.} X. Q. Yin and B. Muckenhoupt,
\newblock  {\em{Weighted inequalities for the maximal geometric mean operator}},
\newblock Proc. Amer. Math. Soc. {\bf{124}}(1996), 75-81.

\bibitem{Zhuang} Y. D. Zhuang,
\newblock  {\em{On inverses of the Holder inequality}},
\newblock J. Math. Anal. Appl. {\bf{161}}(1991), 566-575.

\bibitem{Zuo H. L. and Liu P. D.}
H. L. Zuo and P. D. Liu,
\newblock {\em The minimal operator and weighted inequalities for martingales},
\newblock Acta Math. Scientia. {\bf{26}}(2006), 31-40.
\end{thebibliography}
\end{document}